\documentclass[preprint,sort&compress,a4paper]{elsarticle}
\usepackage{latexsym,a4}

\usepackage{epsf}
\usepackage{color}
\usepackage{epsfig}

\usepackage{graphicx}
\usepackage{xcolor}
\usepackage{amsmath}
\usepackage{latexsym}
\usepackage{amssymb}
\usepackage{amsfonts}

\usepackage{mathdots}
\usepackage{array}

\usepackage{a4}
\usepackage{longtable}

\usepackage{graphicx}
\usepackage{mathdots}
\usepackage{array}

\RequirePackage{hyperref}
\RequirePackage{hyperref}
		\usepackage{pgfplots}
		\usetikzlibrary{plotmarks}
		\usetikzlibrary{external}
		\pgfplotsset{compat=1.5}
		\newlength\figurewidth 
		\newlength\figureheight
		\newlength\markersize
		\tikzset{external/force remake=false} 
\makeatletter
\def\Ddots{\mathinner{\mkern1mu\raise\p@
\vbox{\kern7\p@\hbox{.}}\mkern2mu
\raise3\p@\hbox{.}\mkern2mu\raise5\p@\hbox{.}\mkern1mu}}
\makeatother
\usepackage{booktabs}

\newtheorem{theorem}{Theorem}
\newtheorem{definition}[theorem]{Definition}

\newcommand{\Emin}{E_{\rm min}}
\newcommand{\Emax}{E_{\rm max}}

\newcommand{\e}{\ensuremath{\mathrm{e}}}

\renewcommand{\P}{P}
\newcommand{\R}{\mathbb{R}}
\newcommand{\C}{\mathbb{C}}
\newcommand{\T}{\tau}

\renewcommand{\Re}{\mathrm{Re}}
\renewcommand{\Im}{\mathrm{Im}}
\newcommand{\ii}{\text{i}}

\newcommand{\structurecomment}[1]{}

\def\en{\text{\sc{e}-}}

\begin{document}
\begin{frontmatter}

\title{An efficient algorithm to compute the exponential of skew-Hermitian matrices for the time integration of the Schr\"odinger equation}

\author[rvt]{Philipp Bader}
\ead{bader@uji.es}
\author[els]{Sergio Blanes}
\ead{serblaza@imm.upv.es}
\author[elsa]{Fernando Casas}
\ead{casas@uji.es}
\author[focal,els]{Muaz Seydao\u{g}lu\corref{cor1}}
\ead{m.seydaoglu@alparslan.edu.tr}

\address[rvt]{Departament de Matem\`atiques, Universitat Jaume I, 12071 Castell\'on, Spain.}
\address[els]{Instituto de Matem\'atica Multidisciplinar,
  Universitat Polit\`ecnica de Val\`{e}ncia, E-46022  Valencia, Spain.}
	\address[elsa]{IMAC and Departament de Matem\`atiques, Universitat Jaume I, 12071 Castell\'on, Spain.}
	\cortext[cor1]{Corresponding author}
\address[focal]{Department of Mathematics, Faculty of Art and Science, Mu\c{s} Alparslan University, 49100, Mu\c{s}, Turkey.}

\begin{abstract}
We present a practical algorithm to approximate the exponential of skew-Hermitian matrices up to round-off error
 based on an efficient computation of Chebyshev polynomials of matrices and the corresponding error analysis. It is based on Chebyshev polynomials of degrees 2, 4, 8, 12 and 18  which are computed with only 1, 2, 3, 4 and 5 matrix-matrix products, respectively.
For problems of the form $\exp(-iA)$, with $A$ a real and symmetric matrix, an improved version is presented that computes the sine and cosine of $A$ 
with a reduced computational cost.
 The theoretical analysis, supported by numerical experiments, indicates that the new methods are more efficient than 
schemes based on rational Pad\'e approximants and Taylor polynomials for all
tolerances and time interval lengths. The new procedure is particularly recommended to be used in conjunction with exponential integrators for the 
numerical time integration of the Schr\"odinger equation.  

\end{abstract}


\end{frontmatter}

\section{Introduction}

Given a  skew-Hermitian matrix, $X\in\mathbb{C}^{N\times N}, \ X^H=-X$, we propose in this paper an algorithm to evaluate $\e^X$ 
up to round off accuracy that is more efficient than standard procedures implemented in computing packages for dimensions $N$ 
up to few hundreds or thousands. 

Computing exponentials of skew-Hermitian matrices is very often an intermediate step in the formulation of numerical schemes used
for simulating the evolution of different problems in Quantum Mechanics. Thus, suppose one needs to solve
 numerically the time-dependent Schr\"odinger equation ($\hbar = 1$)
\begin{equation}   \label{Schr0}
  i \frac{\partial}{\partial t} \psi (x,t)  = \hat{H}(t) \psi (x,t), \qquad \psi(x,0)=\psi_0(x).
\end{equation}
Here $\hat{H}(t)$ is in general a time-dependent Hamiltonian operator, $\psi:\mathbb{R}^d\times\mathbb{R}
\longrightarrow \mathbb{C}$ is the wave function representing the
state of the system, and $\psi_0(x)$ is the initial state. 

One possible approach consists in expressing the solution in terms of an orthonormal basis $\{\phi_k(x)\}_{k=1}^{\infty}$ 
that is truncated up to, say, the first $N$ terms. Then, one has
\[
  \psi_0(x)=\sum_{k=1}^N c_k \, \phi_k(x), \qquad \mbox{and} \qquad
  \psi(t,x)=\sum_{k=1}^N c_k(t) \, \phi_k(x), 
\]
where the coefficients $c(t)=(c_1(t),\ldots,c_N(t))^T$ satisfy
\begin{equation} \label{Schr1}
  i \frac{d }{dt} c(t) = {H}(t)\,  c(t), \qquad c(0)=c_{0} \in \mathbb{C}^N,
 \end{equation}
and $H(t)$ is a Hermitian matrix with elements  $H_{\ell m}=\langle \phi_{\ell}|\hat H(t)|\phi_m\rangle, \ \ell,m=1,\ldots,N$. One then subdivides
the time integration interval in a number of subintervals of length $\T$,  and finally computes approximations 
$c_k \simeq c(t_k)$ at times $t_k = k \T$, $k=1,2,3,\ldots$.




Exponential integrators can be used to solve this problem (see \cite{blanes09tme,hochbruck10ei} and references therein) and they require
 the computation at each time step of one or several
matrix exponentials $\e^{-i\, \T\, H_k}$, $k=1,2,\ldots$, where $H_k$ is a Hermitian matrix depending on $H(t)$ at different times.
%
Although efficient algorithms exist to carry out this task by diagonalizing the constant matrix $H_k$, we will show that
 it is indeed possible to compute the exponential 
 in a very efficient way with a different procedure when $\|\T \, H_k\|$ is not too large . This is typically the situation one encounters when exponential integrators are applied to this class of problems \cite{auckenthaler10mea}.
 

The goal of this work is thus
to present an efficient algorithm for computing $\e^{X}$, with $X$ a skew-Hermitian matrix, up to round off accuracy with a minimum number of matrix-matrix products. The algorithm is based on Chebyshev polynomials and an efficient procedure to evaluate polynomials of matrices. If 
$\|X\|$ is large enough,
 this technique can be combined with scaling-and-squaring. Even then, diagonalizing is only superior when a large number of squarings is necessary.

Since the algorithm can also be used to compute $\e^{-i A}$ when $A$ is a Hermitian matrix, just by taking $A = i X$, in the sequel and without
loss of generality we address this problem. 


 Our approach for computing  $\e^{-i \, A}$ is based on approximations  of the form
\begin{equation}   \label{approx.1}
 \e^{-i \, A} \approx \P_m(A),
\end{equation}
where $\P_m(y)$ is a polynomial in $y$ that approximates the exponential $\e^{-i\, y}$. 
 Different choices for such $\P_m(y)$ are available, namely 
truncated Taylor or Chebyshev series expansions 
in an appropriate real interval of $y$. Rational approximations, like Pad\'e approximants, are also a standard technique to compute the exponential
in combination with scaling and squaring 
 \cite{higham05tsa,higham10cma}. In the autonomous case, when $H(t)$ is constant, this is basically equivalent to solve \eqref{Schr1} using a Gauss-Legendre-Runge-Kutta method \cite{dieci94uia} or a Cayley transform \cite{diele98tct}.


Specifically, 
the scaling and squaring technique is based on the property
\begin{equation}  \label{eses}
	\e^{-iA} = \left( \e^{-iA/2^s}  \right)^{2^s} 
	, \quad s\in\mathbb{N}.
\end{equation}
The exponential $\e^{-iA/2^s}$ is then replaced by a polynomial (or rational) approximation $P_m(A/2^s)$. Both
parameters, $s$ and $m$, are determined in such a way that full machine accuracy is achieved with the minimal computational cost. 

An important ingredient in our procedure consists in designing an efficient way to evaluate the approximation $P_m$. In this respect, the technique we propose 
can be considered as a direct descent of the procedure presented in \cite{blanes02hoo} for reducing the number of commutators 
appearing in different exponential integrators. It was later generalized in \cite{bader17aia} to reduce the number of products necessary 
to compute the Taylor polynomials
for approximating the exponential of a generic matrix  (see also \cite{bader19ctm,sastre19btc} for a more detailed treatment).


In fact, the theoretical analysis carried out here and supported by
numerical experiments performed for different Hermitian matrices $A$, indicates that our new schemes are 
more efficient than those based on rational Pad\'e approximants (as used e.g. in {\sc Matlab}) or on Taylor polynomials for all
tolerances. 
The algorithm computes the parameter
\[
  \beta = \|A\|_1
\]
as an upper bound to the spectrum of $A$. As an optional choice, the user can provide upper and lower bounds for the eigenvalues of  the matrix $A$, $\Emin$ and $\Emax$, and this allows one to consider a shift for reducing the overall cost. Then, the algorithm
automatically selects the most efficient polynomial approximation 
for a prescribed error tolerance. 

Although the algorithms based on Taylor polynomial approximations and the use of scaling-and-squaring constructed in \cite{bader19ctm,sastre19btc} can of course be applied
also here, it turns out that in the particular case of skew-Hermitian matrices (with purely imaginary eigenvalues)  it is more convenient instead to apply 
a similar procedure based on 
Chebyshev polynomials. Here only polynomials of degree $m=2,4,8,12$ and $18$ are considered, since  the number of matrix-matrix
products is minimized in those particular cases. Although higher degrees could in principle be taken, it turns out that applying the scaling-and-squaring technique to
lower degree polynomials renders a similar or higher performance.

In many cases, when solving different quantum mechanical or quantum
control problems \cite{auckenthaler10mea} one ends up with a  real and symmetric matrix, $A^T=A\in\R^{N\times N}$, so that
\[
   \e^{-i A} = \cos( A) -i\sin( A),
\]
and we also provide an algorithm for computing $ \cos( A)$ and $\sin( A)$ simultaneously only involving products of real symmetric matrices.
This new algorithm is more efficient than the approach \eqref{approx.1} since that scheme usually requires products of complex matrices,
and other existing algorithms for the simultaneous computation of the matrix sine and cosine \cite{almohy15naf,seydaoglu21ctm}.
The squaring \eqref{eses} 
(also involving products of complex matrices) is then replaced by the double angle formulae 
\[
\cos(2A)=2\cos^2(A)-I=I-2\sin^2(A), \qquad \sin(2A)=2\sin(A)\cos(A),
\]
so that only two products of real symmetric matrices per squaring are required.

In \cite{blanes15aea} an algorithm for approximating $\e^{-i A} v$ for any real symmetric matrix $A$ and any complex vector $v$ was proposed. It is based
 on the idea of splitting and only requires matrix-vector products $A v$ in such a way that the real and imaginary parts of $\e^{- i A} v$
  are approximated in a different way, with a considerable saving in the computational cost with respect to the usual Chebyshev approximation. Here, by contrast,
we focus on problems where the actual computation of $\e^{-iA}$ for any Hermitian matrix is required.  

The plan of the paper is the following. In section \ref{sec.2} we analyze the approximation of the exponential by Taylor and Chebyshev polynomials and by Pad\'e approximants as well as their error bounds. In section \ref{sec.3} we obtain explicitly the Chebyshev polynomials of the degree previously chosen and for the parameters that ensure the error bound previously studied, and next we present the algorithms to evaluate these polynomials with a reduced number of products. The algorithm for the case of a real-symmetric matrix $A$  is also considered. 
Section \ref{sec.4} contains
numerical experiments illustrating the performance of the new methods and some future lines of research are enumerated in the final Section
\ref{sec.5}.



\section{Polynomial approximations}
 \label{sec.2}


Assume that $P_m(y)$ is a $m$th degree polynomial (or a rational function) approximating the function $\e^{-i\, y}$. Then,  the error is bounded 
 (in Euclidean norm)  as
\[
    \|P_m(A) - \e^{-i\, A}\|  \leq \max_{j=0,1,\ldots, N-1} |\P_m(E_j) -\e^{-i\, E_j}| 
\]
in terms of the real eigenvalues $E_0,\ldots,E_{N-1}$ of the Hermitian matrix $A$. If
the spectrum $\sigma(A)=\{E_0,\ldots,E_{N-1}\}$ is contained in an interval of the form $[\Emin,\Emax]$, then
\begin{equation*}
    \|P_m(A) - \e^{-i\, A} \| \leq  \sup_{\Emin \leq y \leq \Emax} |\P_m(y) - \e^{-i\, y}|.
\end{equation*}
The quantities $\Emax$ and $\Emin$ can be estimated in different ways depending on the particular problem
(see e.g. \cite{huang05ase}). Once they have been determined, by introducing the quantities
\begin{equation} \label{shifting}
\alpha = \frac{\Emax + \Emin}{2},
\qquad \beta = \frac{\Emax - \Emin}{2}, \quad \mbox{ and } \quad
   \overline{A} = A - \alpha I,
\end{equation}
it is clear that the spectrum of the shifted operator $\overline{A}$ is contained in an interval centered at the origin, namely
$\sigma(\overline{A}) = \{E_0-\alpha,\ldots,E_{N-1}-\alpha\} \subset [-\beta,\beta]$, so that  
 \begin{equation}   \label{td.2Cheb}
    \e^{-i \, A} = \e^{-i\, \alpha} \, \e^{-i\, \,\beta (\overline{A}/\beta)},
\end{equation}
with $\sigma(\overline{A}/\beta) \subset [-1,1]$.

If the bounds  $\Emin$ and $\Emax$ cannot be estimated in a convenient way, one can always take $\beta=\|A\|_1$, so that $\sigma({A})\leq \beta$, and no shift is considered.

In any event, and without loss of generality, our problem consists now in 
approximating $\e^{-i\,A}$
for a Hermitian matrix $A$ with $\sigma(A) \subset [-\beta,\beta]$ by means of $P_m(A)$. In that case,
\begin{equation}
\label{eq:Err1}
  \frac{\|P_m(A) - \e^{-i\,A} \|}{\| \e^{-i\,A} \|}  = \|P_m(A) - \e^{-i\,A} \| \leq  \epsilon_m(\beta),
\end{equation}
where
\begin{equation}
  \label{eq:Err2}
\epsilon_m(\theta) :=  \sup_{-\theta\leq y \leq \theta}|
P_m (y) - \e^{-i \, y}|
\end{equation}
and $\| \e^{-i\,A} \| = 1$.

\subsection{Taylor polynomial approximation}
 \label{sec.2.1}

An upper bound for the error estimate (\ref{eq:Err2}) of the $m$th degree Taylor polynomial
\begin{equation}\label{eq:Taylor}
 P_m^T(y) \equiv \sum_{k=0}^m\frac{(-i)^k}{k!} y^k
\end{equation}
approximating $\e^{-i\, y}$ can be obtained by computing the Lagrange form of the remainder in the Taylor series
expansion:
\[
  |P_m^T(y) - \e^{-i \, y}| = \frac{1}{(m+1)!} | \e^{-i \xi} (-i y)^{m+1}| = \frac{1}{(m+1)!} |y|^{m+1}
\]
for $\xi\in(0,y)$ so that, from eq. (\ref{eq:Err2}),  
\begin{equation}\label{eq:ErrorTaylor}
 \epsilon_m^T(\theta)  := \frac{\theta^{m+1}}{(m+1)!}.
\end{equation}
Therefore, $P_m^T(A)$ is guaranteed to approximate $\e^{-i\, A}$ up to round-off error as long as $\beta  \le  \theta$ with $\theta$ such that
$ \epsilon_m^T(\theta) \le u=2^{-53}$. We collect in Table~\ref{tab.thetaTaylor} the largest $\theta$ verifying this restriction for the values of $m$
considered in this work. As stated before, only polynomials of degree $m \le 18$ will be employed in practice.

\begin{table}\centering\footnotesize
\caption{\label{tab.thetaTaylor} {\small $\theta$ values for Taylor and Chebyshev polynomials of degree $m$ that can be computed with $\pi$ matrix-matrix products to approximate $\e^{-iA}$ with $A$ skew-Hermitian and guaranteeing that $\epsilon_m(\theta) \le u = 2^{-53}$. 
The $\theta$ value in the column $m=15+$ corresponds to the polynomial of degree 16 built in \cite{sastre19btc} that approximates the
Taylor expansion up to order 15  with 4 products (see section~\ref{sec.5} for details).}}



\

\newcolumntype{H}{@{}>{\lrbox0}l<{\endlrbox}}
\newcolumntype{D}{>{$}r<{$}}
\begin{tabular}{DDDDDDD}
\toprule
	 m:	
		&  	2	&		  4 &   8&    12&  15+&      18
	\\ 
	 \pi:	
		&  	1	&		  2 &   3&    4& 4&      5
	\\ 
	\midrule
\mbox{Taylor pol.}:  & 8.73\en6&1.67\en3&0.0699&0.336&   0.709& 1.147 \\
\mbox{Chebyshev pol.}:& 1.38\en5&2.92\en3&0.1295&0.636&     & 2.212
				\\
\bottomrule
\end{tabular}
\end{table}

\textbf{Remark:} Notice that we can write the polynomial function in the exponential form $P_m^T(\theta) = \e^{-i \, (\theta+\Delta \theta)}$ where $\Delta \theta={\cal O}(\theta^{m+1})$, so condition $|P_m^T(\theta) - \e^{-i \, \theta}|\leq 2^{-53}$ for $\theta\sim{\cal O}(1)$ implies that $|\Delta \theta| \sim 2^{-53}$. However, when backward error analysis is considered one looks for the largest value of $\theta$ such that $\frac{\Delta \theta}{\theta}\leq 2^{-53}$ so different values for $\theta$ are obtained (smaller values when $\theta<1$ and greater values when $\theta>1$).


%
\subsection{Chebyshev polynomial approximation}
 \label{sec.2.2}


The $m$th degree truncation of the Chebyshev series expansion of
$\e^{-i\, y}$ in the interval $y \in [-\theta,\theta]$ reads
\begin{equation}\label{eq:Chebyshev}
 P_{m,\theta}^C(y) :=    J_0(\theta) + 2 \sum_{k=1}^{m} (-i)^k J_k(\theta) \, T_k(y/\theta) ,
\end{equation}
in terms of the Bessel function of the first kind $J_k(t)$ \cite[formula 9.1.21]{abramowitz65hom} and the $k$th Chebyshev polynomial $T_k(x)$ 
generated from the recursion \cite[section 3.11]{olver10nho}
\begin{equation}\label{ChebPolyn}
  T_{k+1}(x)=2xT_k(x)-T_{k-1}(x), \qquad k\geq 1
\end{equation}
with $T_0(x)=1$, $T_1(x)=x$. 

At least three estimates  for $\epsilon_m(\theta)$ may be considered when dealing with Chebyshev polynomial approximations. 
According with the analysis in \cite[section III.2.1]{lubich08fqt}, one can take
\begin{equation}\label{eq:ErrorChebyshev1}
 \epsilon_m^{C_1}(\theta)  :=
4 \left( \e^{1-\theta^2/(2m+2)^2} \frac{\theta}{2m+2}
  \right)^{m+1}.
\end{equation}
On the other hand, in \cite[Theorem 8.2]{trefethen13ata} it is shown that 
\begin{equation}\label{eq:ErrorChebyshev2}
\max_{|y|\leq 1}  \left|e^{\theta y}-P_{m,\theta}^C(y) \right| \leq
  \frac{2M}{\rho^{n}(\rho-1)} = \epsilon_m^{C_2}(\theta),
\end{equation}
where $\displaystyle M=\max_{z\in{\cal E}_{\rho}} \big|e^{\theta z}\big|=e^{\frac{\theta}2(\rho+1/\rho)}$, and ${\cal E}_{\rho}$ denotes the Bernstein ellipse in the complex plane \cite[chapter 8]{trefethen13ata}, 
\[
  {\cal E}_{\rho} =\left\{ z\in \C \ \Big| \ z=\frac12 (r+r^{-1}), \ r=\rho \, \e^{i\phi}, \ 
	-\pi\leq \phi \leq \pi	\right\}.
\]
Here, $\rho$ is any positive number with $\rho>1$, and the optimal value that minimizes the right hand side of \eqref{eq:ErrorChebyshev2} has to be computed numerically for each choice of $\theta$. 

Finally, one can also take the tail of the whole Chebyshev series expansion as an upper bound of the error, i.e.,
\begin{equation}\label{eq:ErrorChebyshev3}
  \|P_{m,\theta}^C(A) - \e^{-i\, A}\| \leq
		\left\| \sum_{k=m+1}^{\infty}  2 (-i)^k J_k(\theta) \, T_k(y/\theta)  \right\| \leq
		\sum_{k=m+1}^{\infty} 2 \, | J_k(\theta) |   \equiv \epsilon_m^{C_3}(\theta). 
\end{equation}
We have evaluated the three estimates (\ref{eq:ErrorChebyshev1})-(\ref{eq:ErrorChebyshev3}) for the relevant degrees $m$ and compared with the 
observed behaviour of the corresponding polynomials. From these computations we conclude that the bound 
 \eqref{eq:ErrorChebyshev3} exhibits the sharpest result, i.e., larger values of $\theta$ for all $m$ considered. Thus, in particular, for $m=18$ bound 
(\ref{eq:ErrorChebyshev1}) leads to $\theta = 1.8843$, bound (\ref{eq:ErrorChebyshev2}) gives $\theta = 1.939$, whereas bound (\ref{eq:ErrorChebyshev3})
provides the largest value $\theta = 2.212$. The corresponding values for $\theta$ obtained with (\ref{eq:ErrorChebyshev3})
are also collected in Table~\ref{tab.thetaTaylor}.
Notice that these values are almost twice larger than those associated to Taylor approximations. 

\begin{figure}[!h]
\begin{center}
	\includegraphics[scale=0.6]{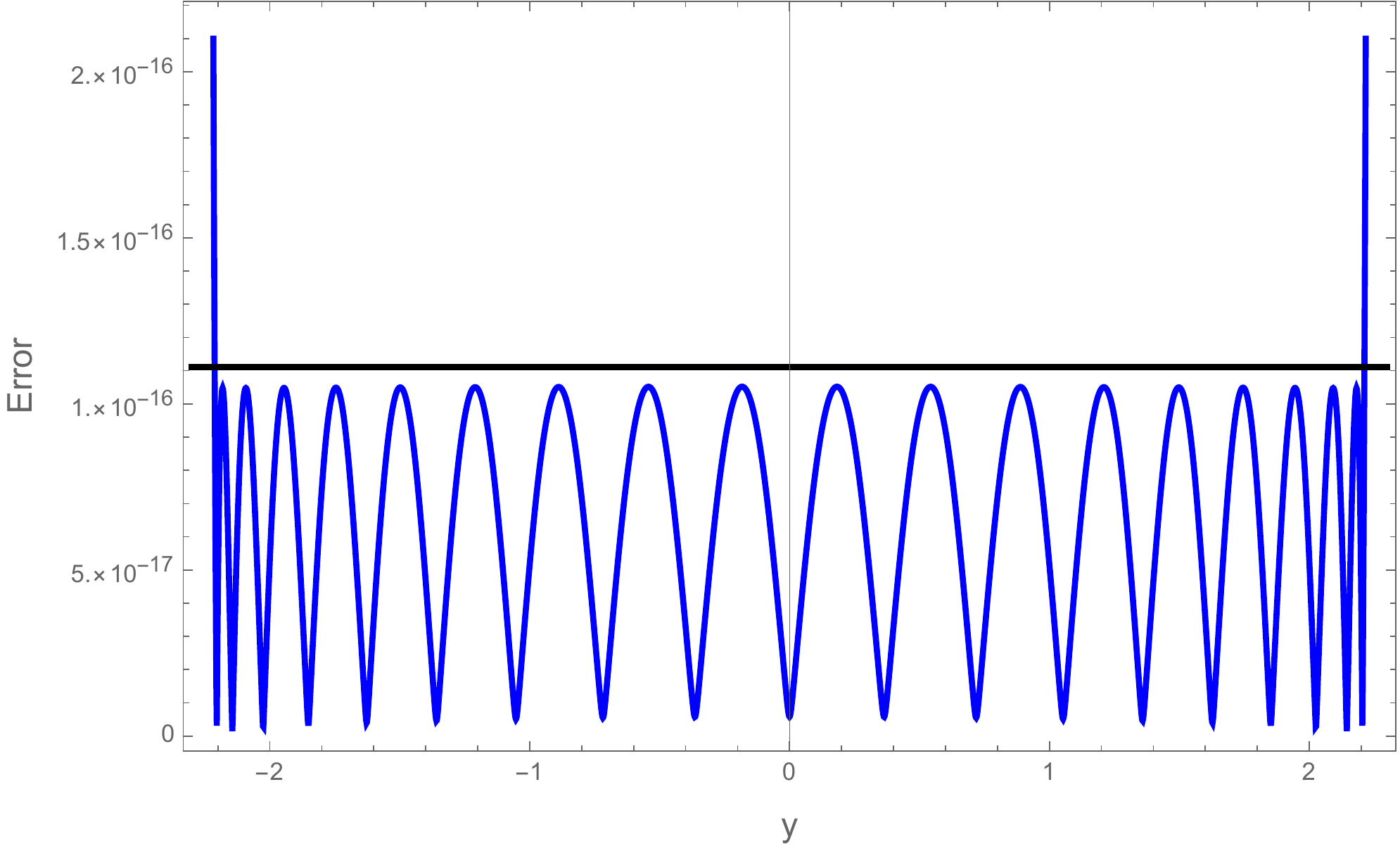}
	\caption{{\small Absolute error $|P_{m,\theta}^C(y) - \e^{-i y}|$ for $(m=18, \theta = 2.212)$  (blue) and the 
	value of $u \approx 1.11\en16$ (black). The error is always smaller than $u$ for $y \in [-\theta,\theta]$.}}	\label{figu1}	 	
\end{center}
\end{figure}

In practice, we have constructed the Chebyshev polynomial approximations for each pair $(m, \theta)$ specified in Table~\ref{tab.thetaTaylor} as
in \cite{gil07nms} (which is in fact equivalent to Eq. \eqref{eq:Chebyshev})
\begin{equation}  \label{ourway1}
  P_{m,\theta}^C(y) = \frac{1}{2} c_0 + \sum_{k=1}^m c_k T_k(y/\theta),
\end{equation}
with
\begin{equation}   \label{ourway2}
   c_k = \frac{2}{\pi} \int_{-1}^1  \frac{\e^{-i \theta y} \, T_k(y)}{\sqrt{1- y^2}} dx
\end{equation}   
and all the calculations have been carried out with 30 digits of accuracy. In Figure \ref{figu1} we show both the absolute error $|P_{m,\theta}^C(y) - \e^{-i y}|$ for 
$(m=18, \theta = 2.212)$ and the value of $u \approx 1.11\en16$. Notice how the error is always smaller than $u$ for the whole interval $y \in [-\theta,\theta]$.


\subsection{Pad\'e approximations}
 \label{sec.2.1b}

Most popular computing packages such as
{\sc Matlab} (\texttt{expm}) and {\it Mathematica} (\texttt{MatrixExp}) use Pad\'e approximants (in combination with scaling-and-squaring) to 
compute numerically the exponential of a generic matrix \cite{higham05tsa,higham10cma}. 

Diagonal {$[m/m]$} Pad\'e approximants are of the form 
\begin{equation}\label{eq.24}
	{	r_{m}(-i\, A )=p_{m}(-iA) \big[ p_{m}(iA ) \big]^{-1} , }
\end{equation}
where 
\begin{equation}  \label{eq.1.2b}
   p_{m}(x)= \sum_{j=0}^m \frac{(2m - j)! m!}{(2m)! (m-j)!} \frac{x^j}{j!},
\end{equation}
and they verify that $r_{m}(-i\, A) = \e^{-i\, A}+\mathcal{O}(A^{2m+1})$. In practice, the evaluation of $p_{m}(-i\, A )$ and $p_{m}(i\, A )$ is carried out 
so as to keep  the
number of matrix products at a minimum. 
The previous notation $\mathcal{O}(A^{n})$ is defined next, since it will be helpful in the sequel.
\begin{definition}
We say that a given function $f(A)$ of the matrix $A$ satisfies $f(A) = \mathcal{O}(A^{n})$ if it can be written as a convergent Taylor expansion,
$f(A) = \sum_{k=n}^{\infty} c_k A^k$, for $\|A\| < \alpha$, with $\alpha$ a positive constant.
\end{definition}
For skew-Hermitian matrices, we can use, instead of the generic backward error bounds obtained e.g. in \cite{higham08fom}, an error estimate of the form
(\ref{eq:Err1}) with $\epsilon_m(\theta)$ in (\ref{eq:Err2}) replaced by its upper bound:
\begin{equation}\label{eq:ErrorPade}
  \|r_m(-i\, A) - \e^{-i\, A}\| \leq
		\left| \sum_{k=2m+1}^{\infty} d_k \theta^k \right| \leq
		\sum_{k=2m+1}^{\infty} |d_k| \, \theta^k  \equiv \epsilon_m^P(\theta). 
\end{equation}
In practice, for a given $m$, we have computed  $s \equiv \sum_{k=2m+1}^{2000} |d_k| y^k$ and determined the largest  $y$ for which $s \le u=2^{-53}$.
This value is taken then as the bound $\theta$. The values for $\theta$ are collected in Table~\ref{tab.thetaPade} for those $m$ for which the diagonal
Pad\'e approximant can be computed with the minimum number of products. 
The function  \verb"expm"  in {\sc Matlab} uses the corresponding bound $\theta$ obtained from relative
backward error with a cost of 2,3,4,5 and 6 products, respectively, in addition to one matrix inverse (we take the cost of one inverse as 4/3 products\footnote{For
a $N \times N$ matrix, it requires one $LU$ factorization at the cost of $1/3$ products plus $N$ solutions of upper and lower triangular systems by
forward and backward substitution at the cost of one product.}). In order to compare with our methods under the same conditions, we have used the  function  \verb"expm" from {\sc Matlab} but taking the $\theta$ values from Table~\ref{tab.thetaPade}.
One should notice that the corresponding backward error bounds  are
smaller up to $m=7$.

\begin{table}\centering\footnotesize
\caption{\label{tab.thetaPade} \small Values of $\theta$ for diagonal Pad\'e approximants of the highest order $2m$ that are computed with $\pi$ products (we take the computation of the inverse of a matris as $4/3$ products).
}
\newcolumntype{H}{@{}>{\lrbox0}l<{\endlrbox}}
\newcolumntype{D}{>{$}r<{$}}
\begin{tabular}{DDDHDHDHDHHHD}
\toprule
	m:	
	&		2
	&		  3 &  		  4&		   5&  		  6&       7&   8&   9&  10&  11&  12&  13 
	\\	
	\pi:	
			&		2+\frac13 
	&		  3+\frac13 &  		  4&		   4+\frac13&  		  6&       5+\frac13&   8&   6+\frac13&  10&  11&  12&  7+\frac13 
	\\ \midrule
\mbox{Pad\'e} &
						2.4007\en3&	
				2.715\en2&1.108\en1 & 2.803\en1  & 0.5443  & 0.8983 & 1.331  &  1.833   & 2.391 & 2.996 & 3.640 & 4.316 
				\\
\bottomrule
\end{tabular}
\end{table}

\section{Evaluating Chebyshev polynomial approximations with a reduced number of products}
\label{sec.3}


Our next goal is to reproduce the Chebyshev polynomial approximations considered in section \ref{sec.2.2} with a reduced number of matrix products
in comparison with the  \emph{de facto} 
standard Paterson--Stockmeyer method for polynomial evaluation. Since the technique has been already explained in detail in
the context of Taylor polynomials approximating the exponential of a generic matrix in \cite{bader19ctm} (see also \cite{sastre18eeo} for a closely related procedure), here we only collect its most salient
features and refer to \cite{bader19ctm} for a comprehensive treatment.

Essentially, the idea is a 
modification of a procedure designed in \cite{blanes02hoo} to reduce the number of commutators appearing in exponential integrators, and consists in taking
a sequence of products of the form
\begin{align}
A_0 &:= I, \qquad A_1 := A \nonumber  \\
A_2 &:= z_0A_0+z_1A_1+\big(x_1A_0+x_2A_1\big)\big(x_3 A_0+x_4 A_1\big) \label{New_algorithm}  \\
A_4 &:= z_2A_0+z_3A_1+z_4A_2+\big(x_5A_0+x_6 A_1+x_7A_2\big)\big(x_8A_0+x_9A_1+x_{10}A_2\big) 
\nonumber \\
A_8 &:= z_5A_0+\cdots+z_8A_4+\big(x_{11}A_0+\cdots+x_{14}A_4\big)\big(x_{15}A_0+\cdots+x_{18}A_4\big),\nonumber \\
 & \vdots  \nonumber
\end{align}
to rewrite any polynomial $P_m(A)$ as $P_m(A) = \sum_{k \ge0} \alpha_k A_k$. Proceeding in this way there might be both
redundancies in the coefficients (for instance, it suffices to take $A_2 = A_1A_1$ since any polynomial of degree two can be writen in terms of $A_0,A_1$ and $A_2$) and also not enough parameters to reproduce some powers in $A$ (e.g. to compute $P_7(A)=A^7$). For this reason, one includes 
new terms of the form, say,
\[
  (w_0A_0+w_1A_1)(w_3A_0+w_4A_1+w_5A_2),
\]
in the procedure for computing $A_k$, {$k > 2$}, so that one has additional parameters. The price to be paid is of course that it is necessary to evaluate
some extra products.

Concerning the particular class of polynomials and degrees we are interested in, $P_m(A)$ with $m=2,4$ can be obtained with just 1 and 2 matrix products,
in a similar way as the Paterson--Stockmeyer technique.

\paragraph{Degree $m=2$} The quadratic Chebyshev polynomial with  $\theta=1.38\en5$ can be trivially computed with one product, and is given by
\[
  A_2 = A^2, \qquad\qquad P_{2,\theta}^C(A)  = \alpha_0 I + \alpha_1 A + \alpha_2 A_2,
\]
with
\[
\begin{aligned}
  \alpha_0 & = 0.9999999999999999999998, \qquad
  \alpha_1  =-0.9999999999761950000001 \, i\\
  \alpha_2 & =-0.4999999999920650000000. 
\end{aligned}
\]
\paragraph{Degree $m=4$} The Chebyshev polynomial of degree four with  $\theta=2.92\en3$ can be computed with two products as follows:
\[
\begin{aligned}
  & A_2  = A^2, \qquad\qquad   A_4 = A_2 (x_1 A + x_2 A_2), \\
  & P_{4,\theta}^C(A) = \alpha_0 I + \alpha_1 A + \alpha_2 A_2 + A_4,
\end{aligned}
\]
with
\[
\begin{array}{lll}
  \alpha_0  = 0.99999999999999999997, & \;\; & \alpha_1  =-0.99999999999981067844 \, i \\
  \alpha_2  =-0.49999999999994320353, &  & x_1     =  0.16666657785001893215 \, i \\
   x_2     =  0.04166664890333648869. &  &
\end{array}
\]  
Although we report here 20 digits for the coefficients, they can be in fact determined with arbitrary accuracy.

The situation is more involved, however, for higher degrees. We next collect the results for
the Chebyshev polynomial approximations to the exponential of degrees $m = 8, 12$ and $18$. Although more values of $m$
could be considered, it turns out that
these polynomials can be constructed with only 3, 4 and 5 products, respectively.

\paragraph{Degree $m=8$}
As is the case with Taylor polynomials \cite{bader19ctm}, the following sequence allows one to evaluate  $P_8(A) \equiv P_{8,\theta}^C(A)$,
with $\theta = 0.1295$:
\begin{equation}  \label{Algorithm83eA}
\begin{aligned}
   & A_2  = A^2,\qquad\qquad   A_4  = A_2(x_1 A+x_2 A_2),\\
   & A_8  = (x_3 A_2+A_4)(x_4I+x_5A+x_6 A_2+x_7A_4),  \\
   & P_{8,0.1295}^C(A)   = \alpha_0 I + \alpha_1 A +\alpha_{2} A_2+ A_8. 
\end{aligned}  
\end{equation}
Notice that this is a particular example of the sequence \eqref{New_algorithm} with some of the coefficients fixed to zero to avoid
redundancies. The parameters $x_i$, $\alpha_i$ are determined such that $P_{8,0.1295}^C(A)$ agrees with the corresponding expression (\ref{ourway1}).
One has 10 parameters to solve 9 nonlinear equations and this results in two families of solutions depending on a free parameter, $x_1$. All solutions provide the same polynomial (if exact arithmetic is considered), and we have chosen $x_1$ to (approximately) minimize the 1-norm of the vector of parameters. The corresponding 
coefficients in (\ref{Algorithm83eA}) for the Chebyshev polynomial are given by
\[
\begin{array}{lll}
x_1=431/4000,   & \;\; &
x_2=-0.02693906873598870733 \, i , \\
x_3= 0.66321004441662438593 \, i , &  &
x_4= 0.54960853911436015786 \, i, \\
x_5= 0.16200952846773660904,&  &
x_6=-0.01417981805211804396 \, i , \\
x_7=-0.03415953916892111403,&  &
\alpha_0= 0.99999999999999999928, \\
\alpha_1=-0.99999999999999233987 \, i,&  &
\alpha_2=-0.13549409636220703066.
\end{array}
\]

\paragraph{Degree $m=12$}
Here the situation is identical to what happens with Taylor polynomials approximating $\e^X$ for a generic matrix \cite{bader19ctm}: 
although polynomials up to degree 16 could in principle be constructed with 4 products {by applying the sequence \eqref{New_algorithm}}, 
in practice the highest degree we
are able to get is $m=12$ with the following sequence:
\begin{equation}  \label{poly12}
	\begin{array}{lll}
		A_2=A^2, & \;\; & 	A_3=A_2A, \\
		B_1 = a_{0,1}I+a_{1,1}A+a_{2,1}A_2+a_{3,1}A_3, & & 
		B_2 = a_{0,2}I+a_{1,2}A+a_{2,2}A_2+a_{3,2}A_3,\\
		B_3 = a_{0,3}I+a_{1,3}A+a_{2,3}A_2+a_{3,3}A_3, & & 
		B_4 = a_{0,4}I+a_{1,4}A+a_{2,4}A_2+a_{3,4}A_3,\\
		A_6 = B_3 + B_4^2 & & \\
		P_{12,0.636}^C(A)  = B_1 + (B_2 + A_6)A_6. & & 
	\end{array}
	\end{equation}
This ansatz has four families of solutions with three free parameters. A judicious choice leading to a small value for 
 $\sum_{i,j} |a_{i,j}|$ is:
\[
		\begin{array}{lll}
	a_{0,1} =-6.26756985350202252845, &  \;\; &
	a_{1,1} = 2.52179694712098096140  \, i, \\
	a_{2,1} = 0.05786296656487001838, &  &
	a_{3,1} =-0.07766686408071870344 \, i, \\
	a_{0,2} = 0,  &  &
	a_{1,2} = 1.41183797496250375498 \, i, \\%
	a_{2,2} = 0,  &   &
	a_{3,2} =-0.00866935318616372016 \, i, \\
	a_{0,3} = 2.69584306915332564689, & &
	a_{1,3} =-1.35910926168869260391  \, i, \\%
	a_{2,3} =-0.09896214548845831754, & &
	a_{3,3} = 0.01596479463299466666  \, i, \\%
	a_{0,4} = 0, &  &
	a_{1,4} = 0.13340427306445612526 \, i, \\
	a_{2,4} = 0.02022602029818310774, &   &
	a_{3,4} =-0.00674638241111650999 \, i.
	\end{array}
	\]

\paragraph{Degree $m=18$}
We have been able to write the Chebyshev polynomial approximation of degree $m=18$ with 5 products.  This is done by expressing $P_{18,2.212}^C(A)$ as the product of two polynomials of degree 9, 
that are further decomposed into polynomials of lower degree. The polynomial is evaluated through the following sequence:
\begin{equation}   \label{poly.18}
\begin{aligned}
& A_2 = A^2, \qquad\qquad A_3 = A_2A, \qquad\qquad A_6 = A_3^2,\\
& B_{1} = a_{0,1}I + a_{1,1}A + a_{2,1}A_2 + a_{3,1}A_3,\\
& B_{2} = b_{0,1} I + b_{1,1}A + b_{2,1}A_2 + b_{3,1}A_3 + b_{6,1}A_6, \\
& B_{3} = b_{0,2} I + b_{1,2}A + b_{2,2}A_2 + b_{3,2}A_3 + b_{6,2}A_6,\\
& B_{4} = b_{0,3} I + b_{1,3}A + b_{2,3}A_2 + b_{3,3}A_3 + b_{6,3}A_6,\\
& B_{5} = b_{0,4} I+  b_{1,4}A + b_{2,4}A_2 + b_{3,4}A_3 + b_{6,4}A_6,\\
& A_{9}  = B_{1}B_{5} + B_{4},\\
& P_{18,2.212}^C(A) = B_{2} + (B_{3} + A_9)A_9,
\end{aligned}								
\end{equation}
with coefficients 
\[
\begin{array}{lll}
  a_{0,1} = 0,  &  \;\;  &  a_{1,1} =  3/25,\\
  a_{2,1} = -0.00877476096879703859 \, i,  &  &
  a_{3,1} = -0.00097848453523780954, \\
  b_{0,1} =  0, &  &
  b_{1,1} = -0.66040840760771318751 \, i, \\
  b_{2,1} = -1.09302278471564897987, & &
  b_{3,1} =  0.25377155817710873323 \, i, \\
  b_{6,1} =  0.00054374267434731225, & &
  b_{0,2} = -2.58175430371188142440, \\
  b_{1,2} = -1.73033278310812419209 \, i, & &
  b_{2,2} = -0.07673476833423340755, \\ 
  b_{3,2} = -0.00261502969893897079 \, i, &  &
  b_{6,2} = -0.00003400011993049304, \\
  b_{0,3} =  2.92377758396553673559, & &
  b_{1,3} =  1.44513300347488268510 \, i, \\
  b_{2,3} =  0.12408183566550450221, &  &
  b_{3,3} = -0.01957157093642723948 \, i, \\
  b_{6,3} =  0.00002425253007433925, &  &
  b_{0,4} = 0, \\
  b_{1,4} = 0, &  &
  b_{2,4} = -0.123953695858283131480 \, i, \\
  b_{3,4} = -0.011202694841085592373,  &   &
  b_{6,4} = -0.000012367240538259896 \, i. 
\end{array}
\]

\subsection{The case of real symmetric matrices}

In the particular case when $A$ is a real symmetric matrix, we can write
\begin{equation} \label{eq:unitary2}
  \e^{-i A} = \cos( A) -i\sin( A)
\end{equation}
and it is possible to construct algorithms
for approximating the real symmetric matrices $ \cos(A)$ and $\sin(A)$ simultaneously by means of a reduced number of 
products of real symmetric matrices,  as shown in \cite{almohy15naf,seydaoglu21ctm}.

The polynomial (\ref{ourway1}) can be decomposed into real and imaginary parts,
\[
  P^C_{m,\theta}(A) = c^C_{m,\theta}(A) - i s^C_{m,\theta}(A),
\]
with
%
\[
  c^C_{m,\theta}(A)=\Re(P^C_{m,\theta}(A)), \qquad\qquad
  s^C_{m,\theta}(A)=-\Im(P^C_{m,\theta}(A))
\]
and the goal is to compute exactly $c^C_{m,\theta}(A)$ with a reduced number of products. Then, by using all computations carried out in this process, one also obtain approximations $S^C_{m,\theta}(A)$  to the imaginary part in such a way that {it is a polynomial of degree $k>m$ such that} $S^C_{m,\theta}(A)=s^C_{m,\theta}(A)+{\cal O}(A^{m+1})$.
Taking into account that
\begin{eqnarray}\label{eq:ErrorChebyshevSym}
  \|\Big(c^C_{m,\theta}(A)+iS^C_{m,\theta}(A)\Big) - \e^{-i\, A}\| 	
	& \leq &
  \|P^C_{m,\theta}(A) - \e^{-i\, A}\| +
  \|s^C_{m,\theta}(A)-S^C_{m,\theta}(A)\| \nonumber   \\
	& \leq &
 \epsilon_m^{C_3}(\theta)+ 
  \|s^C_{m,\theta}(A)-S^C_{m,\theta}(A)\| 
\end{eqnarray}
we need to check if 
$\|s^C_{m,\theta}(A)-S^C_{m,\theta}(A)\|\leq\sum_{\ell=m+1}^k |c_{\ell}|\, |\theta|^{\ell}\leq 2^{-53}$
for these values of $m$ and $\theta$. If this is not the case, one has to find the maximum value of $\vartheta$ such that $\|s^C_{m,\theta}(\vartheta)-S^C_{m,\theta}(\vartheta)\| \le 2^{-53}$ and to take this value as the value for $\theta$, i.e. the largest value of $\|A\|$ that guarantees an error smaller that roundoff. 
If the value obtained for $\vartheta$ is considerably smaller than the value of $\theta$ obtained for the cosine function,
we will look for a new and more accurate approximation to $s^C_{m,\theta}(A)$ by taking e.g. one extra product in the numerical scheme. 



For example, we can compute simultaneously $c^C_{4,\theta}(A)$ and $s^C_{4,\theta}(A)$, i.e. $P^C_{4,\theta}(A)$, with three products of real symmetric matrices. However, with the same number of products we can also compute $c^C_{5,\theta}(A), \ s^C_{5,\theta}(A)$, i.e. $P^C_{5,\theta}(A)$ that has a larger value of $\theta$, so that we only consider this last case which is evaluated as follows (in this case $S^C_{m,\theta}(A)=s^C_{m,\theta}(A)$).
 

\paragraph{Degree $m=5$}
The polynomial $c_{5,\theta}^C(A)$, 
$\theta = 1.17$E-2, 
is computed with 2 products by taking $B = A^2$ as
\[
   c_{5,\theta}^C(A)   = \alpha_0 I + \alpha_1 B +\alpha_{2} B^2,
\]
with
\[
\begin{array}{lll}
  \alpha_0  = 0.99999999999999988866, & \;\; &
  \alpha_1  =-0.49999999998536031183 \\
  \alpha_2  = 0.04166638147997997916, & & 
\end{array}
\]
whereas for evaluating $s_{5,\theta}^C(A)$  only one additional product is required:
\[
    s_{5,\theta}^C(A)  = A( z_0 I + z_1 B + z_2 B^2)
\]
with
\[
\begin{array}{lll}
  z_0 = 0.99999999999999994433 , & \;\; &
  z_1 =-0.16666666666341340086 \\
  z_2 = 0.00833328580219952161. & & 
\end{array}
\]  

\paragraph{Degree $m=8$}
The polynomial $c_{8,\theta}^C(A)$, $\theta = 0.1295$, 
is computed with 3 products as:
\[
\begin{aligned}
  & B = A^2 ,\qquad\qquad B_2  = B^2,\qquad
  B_4  = B_2(x_1 B+x_2 B_2),\\
   & c_{8,\theta}^C(A)   = \alpha_0 I + \alpha_1 B +\alpha_{2} B_2+ B_4,
\end{aligned}
\]
with
\[
\begin{array}{lll}
  \alpha_0  = 0.99999999999999999928, & \;\; &
  \alpha_1  =-0.49999999999999787210, \\
  \alpha_2  = 0.04166666666565156615, & & 
   x_1     = -0.00138888871939942118,  \\
   x_2     =  0.00002479003614491668, & & 
\end{array}
\]  
and $s_{8,\theta}^C(A)$ is approximated with error $\mathcal{O}(A^9)$ with one additional product by
\[
  S_{8,\theta}^C(A) =   A (z_0 I + z_1 B + z_2 B_2 + z_3 c_{8,\theta}^C(A) )
\]
 with
\[
\begin{array}{lll}
   z_0     =   0.85721768947064012466, & \;\; &
   z_1     =  -0.09527551139590047256, \\
   z_2     =   0.00238406908730568850, & & 
   z_3     =   0.14278231052935221530.
\end{array}
\]
Notice that the condition
\[
  |s^C_{8,\theta}(\vartheta)-S^C_{8,\theta}(\vartheta)|= x_2 \, z_3 \, \vartheta^9 \le 2^{-53}
\]
is satisfied only for $\vartheta \le 0.06807$. This is a significant reduction with respect to $\theta$ and for this reason we look for an approximation which involves one extra product. With five products, however, it is possible to exactly compute the polynomials for $m=9$, with $\theta = 0.2143$. 

\paragraph{Degree $m=9$}
The polynomial $c_{9,\theta}^C(A)$, 
$\theta = 0.2143$, 
is computed with 4 products in the same way:
\[
\begin{aligned}
  & B = A^2,\qquad\quad   B_2  = B^2,\qquad\quad   B_3  = B_2B,\qquad\quad 
  B_4  = B_3B,\\
   & c_{9,\theta}^C(A)   = \alpha_0 I + \alpha_1 B +\alpha_{2} B_2+\alpha_{3} B_3+\alpha_{4} B_4,
\end{aligned}
\]
with
\[
\begin{array}{lll}
  \alpha_0  = 0.99999999999999989168, & \;\; &
  \alpha_1  =-0.49999999999988173685,   \\
  \alpha_2  = 0.04166666664600636231, & & 
  \alpha_3  =-0.00138888762558264513,   \\
  \alpha_4  = 0.00002477005498155486,  & & 
\end{array}
\]  
and $s_{9,\theta}^C(A)$ can be computed with one additional product:
\[
\begin{aligned}
   &  s_{9,\theta}^C(A) =   A (z_0 I + z_1 B + z_2 B_2 + z_3 B_3 + z_4 B_4 ),
\end{aligned}
\]
 with
\[
\begin{array}{lll}
   z_0     = -0.999999999999999945837  , & \;\; &
   z_1     =  0.166666666666643012068, \\
   z_2     = -0.008333333330440664914  , & & 
   z_3    =   0.000198412554024823435 ,\\
   z_4     = -2.75257852630876250884\cdot 10^{-6} .  & & 
\end{array}
\]

\paragraph{Degree $m=16$}
We can compute $c_{16,\theta}^C(A)$, $\theta = 1.5867$, with only four products as follows. 
We first take $B=A^2$, so that $c_{16,\theta}^C$ is indeed a polynomial of degree eight in $B$, that can be computed with only three products in a similar way to
$P_{8,0.1295}^C(A)$ with the sequence
\begin{equation}  \label{Algorithm83eAc}
\begin{aligned}
   &  B=A^2,\qquad\qquad B_2  = B^2,\qquad\qquad
   B_4  = B_2(x_1 B+x_2 B_2),\\
   & B_8  = (x_3 B_2+B_4)(x_4I+x_5B+x_6 B_2+x_7B_4),  \\
   & c_{16,\theta}^C(A)   = \alpha_0 I + \alpha_1 B +\alpha_{2} B_2+ B_8,
\end{aligned}  
\end{equation}
where
\[
\begin{array}{lll}
x_1=1/100, &  \;\;  &
x_2=-0.00008035854055477845 , \\
x_3=-0.10743065643419630630,  &  &  
x_4=-0.12491372919298427513 , \\  
x_5= 0.00130085397953037838 ,&  &
x_6=-0.00001633763177694857 , \\
x_7= 7.13215089463286614820 \cdot 10^{-6}, &  &
\alpha_0= 0.99999999999999999530, \\
\alpha_1=-0.49999999999999969795, &  &
\alpha_2=0.028247102741817734721.
\end{array}
\]
With two extra products we can approximate the matrix $s_{16,\theta}^C(A)$:
\begin{equation}  \label{Algorithm6sA}
\begin{aligned}
   C_{24} & = (z_5 I + z_5 B +z_{6} B_2+ z_{7} B_4+ z_{8} \, c_{16,\theta}(A))B_4,  \\
   S_{16,\theta}^C(A)  & = A\big(z_0 I + z_1 B +z_{2} B_2+ z_{3} B_4+ z_{4} \, c_{16,\theta}(A) + C_{24}\big), 
\end{aligned}  
\end{equation}
with
\[
\begin{array}{lll}
z_0= 33/50, & \;\; &
z_1=  0.00333333333335438849, \\
z_2= -0.00583333333345309522, &  & 
z_3=  0.02773310749258735833,\\
z_4=  0.33999999999999886261,&  & 
z_5= -0.00034915267907803119,\\
z_6=  4.19573036995827807213\cdot 10^{-6} ,&  & 
z_7= -2.63931697420854364428\cdot 10^{-6},\\
z_8 =-3.00240279002259730782\cdot 10^{-6}.  &  & 
\end{array}
\]
In this way $S_{16,\theta}^C(A)$ is a polynomial of degree 25 in $A$ where the condition
\[
  \|s_{16,\theta}^C(\vartheta)-S_{16,\theta}^C(\vartheta)\| \le 2^{-53}
\]
is satisfied for $\vartheta \le 0.7563$. One extra product (7 products in total) suffices to exactly compute $s_{16,\theta}^C(A)$ (and then to keep the value of $\theta$).
We do not show this scheme because, as we will see, with 7 products one can find an improved approximation.

\paragraph{Degree $m=24$}
The same strategy can be applied to the polynomials $c_{24,\theta}^C(A)$ and $s_{24,\theta}^C(A)$, with $\theta = 4.5743$. Thus, 
$c_{24,\theta}^C(A)$ is computed by taking $B=A^2$ and computing the corresponding polynomial of degree 12 with only four additional products as previously:
\begin{equation}  \label{poly12b}
	\begin{array}{lll}
		D=A^2, \qquad\quad D_2=D^2 & \; &	 D_3=D_2D, \\
		B_1 = a_{0,1}I+a_{1,1}D+a_{2,1}D_2+a_{3,1}D_3, & & 
		B_2 = a_{0,2}I+a_{1,2}D+a_{2,2}D_2+a_{3,2}D_3,\\
		B_3 = a_{0,3}I+a_{1,3}D+a_{2,3}D_2+a_{3,3}D_3, & & 
		B_4 = a_{0,4}I+a_{1,4}D+a_{2,4}D_2+a_{3,4}D_3,\\
		D_6 = B_3 + B_4^2  & & \\
		c_{24,\theta}^C(A)  = B_1 + (B_2 + D_6)D_6, & & 
	\end{array}
	\end{equation}
with
\begin{equation}  \label{Cocos24}
		\begin{array}{lll}
	a_{0,1} =  0.39272620931352327385, &  \;\;&  
	a_{1,1} = -0.08760637124112618048, \\
	a_{2,1} =  0.01962064507143601071,  &  &
	a_{3,1} = -0.00013421604022829771, \\
	a_{0,2} =  1/5,&  &
	a_{1,2} = -0.54235659842328961975, \\%
	a_{2,2} =  679/100000,&  &
	a_{3,2} = -0.00002902999756981724, \\
	a_{0,3} =  0.68566773555140770915, &  &
	a_{1,3} = -0.02578520551577453856, \\%
	a_{2,3} =  0.00019815665089300452,&  &
	a_{3,3} = -1.10083330495602029332\cdot10^{-6}, \\%
	a_{0,4} =  0,&  &
	a_{1,4} = -0.03931944346958836562, \\
	a_{2,4} =  0.00017839382197658767,&  &
	a_{3,4} = -1.06908694221941432625\cdot10^{-6},
	\end{array}
	\end{equation}
whereas with two extra products the following approximation to $s_{24,\theta}^C(A)$ is obtained 
\begin{equation}  \label{AlgSin}
\begin{aligned}
  & C_{48}  = (z_6 I + z_7 D +z_{8} D_2+ z_9 D_3 + z_{10} D_6+ z_{11} \, c_{24,\theta}(A)) c_{24,\theta}(A),  \\
  & S_{24,\theta}^C(A)   = A\big(z_0 I + z_1 D +z_{2} D_2+ z_{3} D_3+ z_{4} D_6 + z_5 \, c_{24,\theta}(A) + C_{48}\big)
\end{aligned}  
\end{equation}
with
\[
\begin{array}{lll}
z_0= -0.01238438326981811663, &  \; \; & 
z_1= -0.06180067679127220638, \\
z_2=  0.00046275599640408615, &  & 
z_3= -9.92990416300441584763\cdot10^{-6},\\
z_4=  1.26307934615308708610,&  & 
z_5=  9.10439014880980346565\cdot10^{-15},\\
z_6=  0.14610549096048524519,&  & 
z_7=  0.00087697762149660844,\\
z_8 = 4.12092186281469998191\cdot10^{-6}, &  & 
z_9 = 2.23743615053828476204\cdot10^{-8}, \\
z_{10} = 0.00033015662857238333, &  & 
z_{11} = -2.405371071766852323329\cdot10^{-7}.
\end{array}
\]
The approximation $S_{24,\theta}^C(A)$ given by (\ref{AlgSin}) is a polynomial of degree 48 in $A$ verifying the condition
\[
  \|s_{24,\theta}^C(\vartheta)-S_{24,\theta}^C(\vartheta)\| \le 2^{-53}
\]
for $\vartheta \le 2.1556$, which is smaller than the value of $\theta$ for this case but larger than the value of $\theta$ for $m=16$ that requires the same number of products, and for this reason the previous scheme is not considered in practice. 

One extra product suffices to exactly compute $s_{24,\theta}^C(A)$ (and then to keep the value of $\theta$) as follows
\begin{equation}  \label{AlgSin2}
\begin{aligned}
  & D_{5}  = D_2(z_{11}D_2+z_{12}D_3),  \\
  & C_{24}  = (z_6 I + z_7 D +z_{8} D_2+ z_9 D_3 + D_5+ z_{13} D_6) (D_6 + z_{10} D),  \\
  & s_{24,\theta}^C(A)   = A\big(z_0 I + z_1 D +z_{2} D_2+ z_{3} D_3+ z_{4} D_5 + z_5 \, c_{24,\theta}(A) + C_{24}\big)
\end{aligned}  
\end{equation}
In this case we can solve all the equations (including the corresponding to $A^{25}$). We have now one free parameter and one solution is:
\[
\begin{array}{lll}
  z_0=  2.85247650396873609664, &  \; \; & 
  z_1= -0.23838922984354509797,  \\
  z_2=  0.01254735251131974478, &  & 
  z_3= -0.00003184984233834954,  \\
  z_4= -7.91411934357932811110,&  & 
  z_5= -0.45584956828766694538,  \\
  z_6= -2.34944723110594310069,&  & 
  z_7= -0.34315650534099675485,  \\
  z_8 = 0.00379529409295014610, &  & 
  z_9 =-0.00001509312002244718, \\
z_{10}=-17/1000, &  & 
z_{11}= 7.68145795118100472945 \cdot 10^{-9},  \\
z_{12}=-2.71896175810263278764\cdot 10^{-11} &  & 
z_{13}= 0.45584956828766694538
%
\end{array}
\]


Table \ref{tab.thetaCosSinCheb} collects the values of $\theta$ for the selected approximations to the sine and consine functions and their cost in terms of products of real symmetric matrices.

\begin{table}\centering\footnotesize
\caption{\label{tab.thetaCosSinCheb} {\small $\theta$ values for Chebyshev polynomials of degree $m$ that can be computed with $\pi$ products of symmetric real matrices to simultaneously compute the sine and cosine matrix functions to approximate $e^{-iA}$ with $A$ a real symmetric matrix (in parenthesis it is indicated the cost and maximum value of $\|A\|$ when the Chebyshev polynomial for the sine function is approximated with a higher degree polynomial).}
}

\

\newcolumntype{H}{@{}>{\lrbox0}l<{\endlrbox}}
\newcolumntype{D}{>{$}r<{$}}
\begin{tabular}{DDDDDD}
\toprule
	 m:	
		&  	5	&		  8 &   9&    16&      24
	\\ 
	 \pi:	
		&  	3	&		  4 &   5&    (6)7&      (7)8
	\\ 
	\midrule
\mbox{Chebyshev pol.}:& 1.17\en2&0.068&0.214&(0.7563)1.587& (2.1556)4.574
				\\
\bottomrule
\end{tabular}
\end{table}

\subsection{The algorithm}
\label{flowalgo}

In previous sections
we have computed a number of Chebyshev polynomials of different degrees for some values of $\theta$ that provide errors below roundoff when approximating $\e^{-i y}$ for $y\in[-\theta,\theta]$. These polynomials are computed by applying a particular sequence in order to reduce the number of products. 
To approximate $\e^{-i\, A}$ one has to select the most appropriate polynomial that leads to an error below the prescribed
tolerance at the smallest computational cost.

The user has to provide the matrix  $A$ and, as an optional input, the values for $\Emin$ and $\Emax$. The algorithm  then computes $\beta$ and determines the normalized matrix $A$. If $\Emin$ and $\Emax$ are not given, the algorithm takes  $\beta=\|A\|_1$ as an upper bound to  $|\Emin|$ and $|\Emax|$ and no shift is considered.

Next, the algorithm determines the most efficient method (among the list of available schemes) leading to the desired result: it chooses the cheapest method with error bounds below round off error. 

If none of the methods provides an error below tolerance, then the scaling and squaring technique is used. In that case, the value of $\theta$ for the polynomial of the highest degree 18 (or 24 for the trigonometric matrix functions) is taken to obtain the number of squarings that will be necessary. 

As an illustration, suppose one is interested in computing $\e^{-iA}$, where $A$ is a complex Hermitian matrix such that $\Emin$, $\Emax$  are not known and,
in addition
\begin{enumerate}
	\item $\|A\|_1=8$. (i) With Pad\'e one checks that $\|A/2\|_1=4<4.316$ and the exponential is computed with one scaling and the approximant with $m=13$ that involves 6 products and one inverse ($8+1/3$ products in total). (ii) With Taylor we have $\|A/2^3\|_1=1<1.1468$, the exponential is computed with three scalings and the polynomial with $m=18$, requiring 5 products (for a total of $8$ products). Finally, (iii)  with Chebyshev, since $\|A/2^2\|_1=2<2.212$, the exponential is computed with two scalings and 5 products (for a total of $7$ products).
	\item $\|A\|_1=0.1$. (i) With Pad\'e, the exponential is computed with 3 products and one inverse; (ii) Taylor requires 4 products, and (iii)  Chebyshev needs
	only 3 products.
	\item $\|A\|_1=0.0025$. (i) With Pad\'e, the exponential is computed with 2 products and one inverse; (ii) Taylor requires 3 products, and (iii) Chebyshev needs 2 products.
\end{enumerate}
Notice that, whereas the reduction in computation is roughly the same in all cases, the relative saving increases as the norm of the matrix is smaller.

This strategy has been implemented as a {\sc Matlab} code which is freely available for download at the website \cite{website21}, together with some notes and examples illustrating the whole procedure.

\section{Numerical examples}
\label{sec.4}

In this section we report on two numerical experiments carried out by applying the previous algorithm based on Chebyshev polynomials. We also
compare their main features with Taylor polynomials and Pad\'e approximants.


\paragraph{Example 1: A high dimensional Rosen--Zener model.} 

This is a generalization of the well known Rosen--Zener model for a quantum  system of two levels \cite{rosen32dsg} which is closely related
to the problem analyzed in \cite{kyoseva07pro}. The corresponding Schr{\"o}dinger equation  (\ref{Schr1}) for the evolution operator (in the interaction picture) is
\begin{subequations}
\label{eq:HIs}
\begin{equation}
\begin{cases}
U'(t) = 
- \, \ii \, H(t) \, U(t)\,, \qquad t \in (t_0, t_f)\,, \\
U(t_0) = I,\,
\end{cases}
\end{equation}
where the time-dependent Hamiltonian reads, after normalization,
\begin{equation}
\begin{gathered}
H(t) = f_1(t) \, \sigma_1 \otimes I  + f_2(t) \, \sigma_2 \otimes R  \in \C^{d \times d}\,, \qquad d = 2 \, k\,, \\
\end{gathered}
\end{equation}
with the identity matrix $I \in \R^{k \times k}$, Pauli matrices 
\begin{equation}
\sigma_1 = \begin{pmatrix} 0 & \;\; \qquad & 1 \\ 1 & & 0 \end{pmatrix}\,, \qquad 
\sigma_2 = \begin{pmatrix} 0 & \quad & - \, i \\ i & & 0 \end{pmatrix}\,,
\end{equation}
and
\begin{equation}
\begin{gathered}
R = \mbox{tridiag}\big(1, 0, 1\big) \in \R^{k \times k}\,.
\end{gathered}
\end{equation}
We take in particular 
\begin{equation}
f_1(t) = V_0 \, \cos(\omega \, t) \, \big(\cosh\big(\tfrac{t}{T_0}\big)\big)^{-1}\,, \qquad 
f_2(t) = - \, V_0 \, \sin(\omega \, t) \big(\cosh\big(\tfrac{t}{T_0}\big)\big)^{-1}\,, 
\end{equation}
\end{subequations}
with $\omega = 5$, $d = 20$, $V_0 = 2$.  We then integrate from $t_0 = -4 T_0$ until the final time $t_f = 4 \, T_0$ for $T_0=1$ and determine numerical approximations, $U_{\text{app}}(t_f,t_0)$ at $t = t_f$ for different time step sizes $\tau = \frac{t_f - t_0}{M}$; 
 (a reference solution~$U_{\text{ref}}(t_f,t_0)$ is computed numerically to high accuracy).


%

\begin{figure}[!ht] 
\centering
  \includegraphics[width=.85\textwidth]{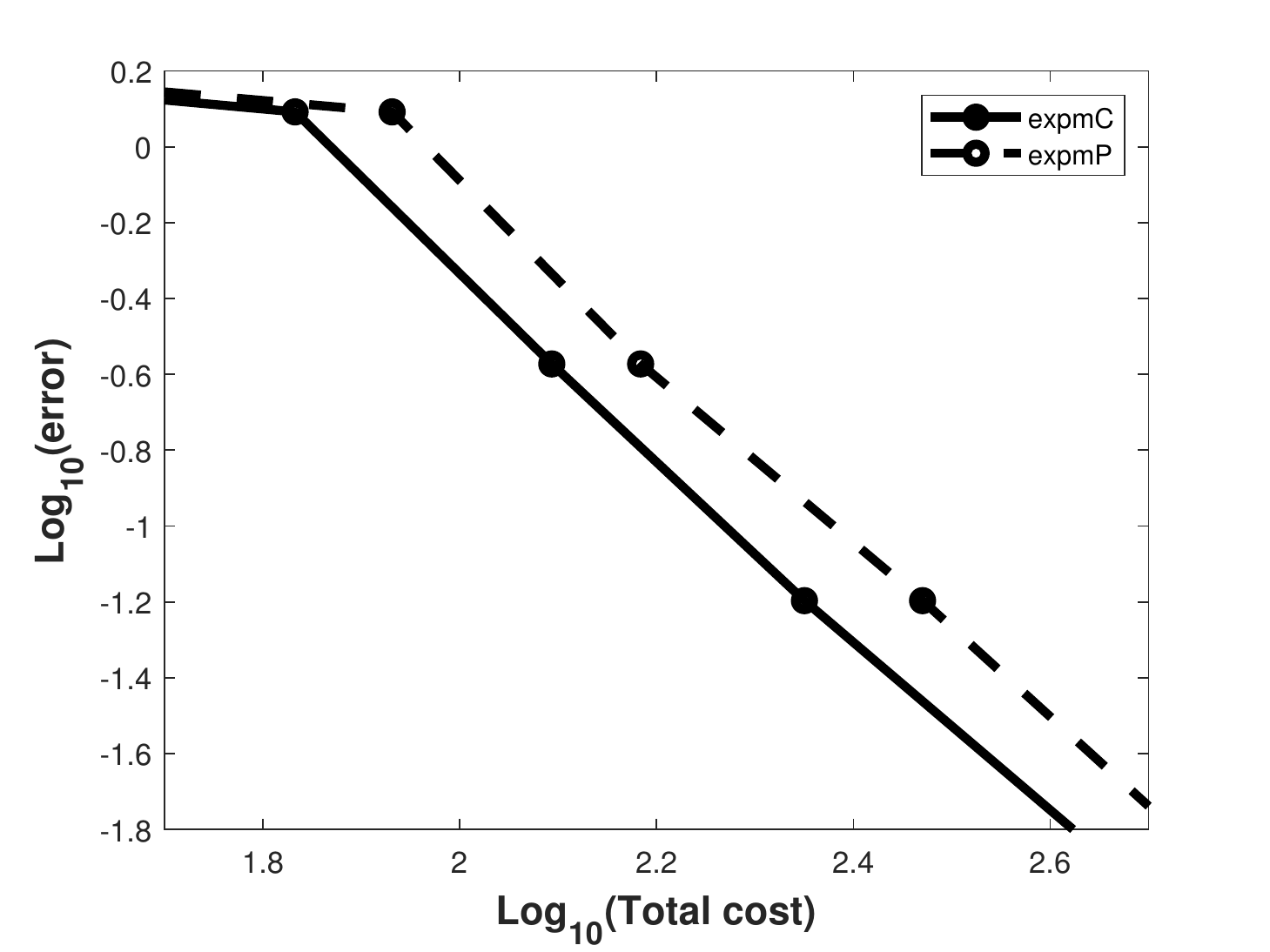}
  \includegraphics[width=.85\textwidth]{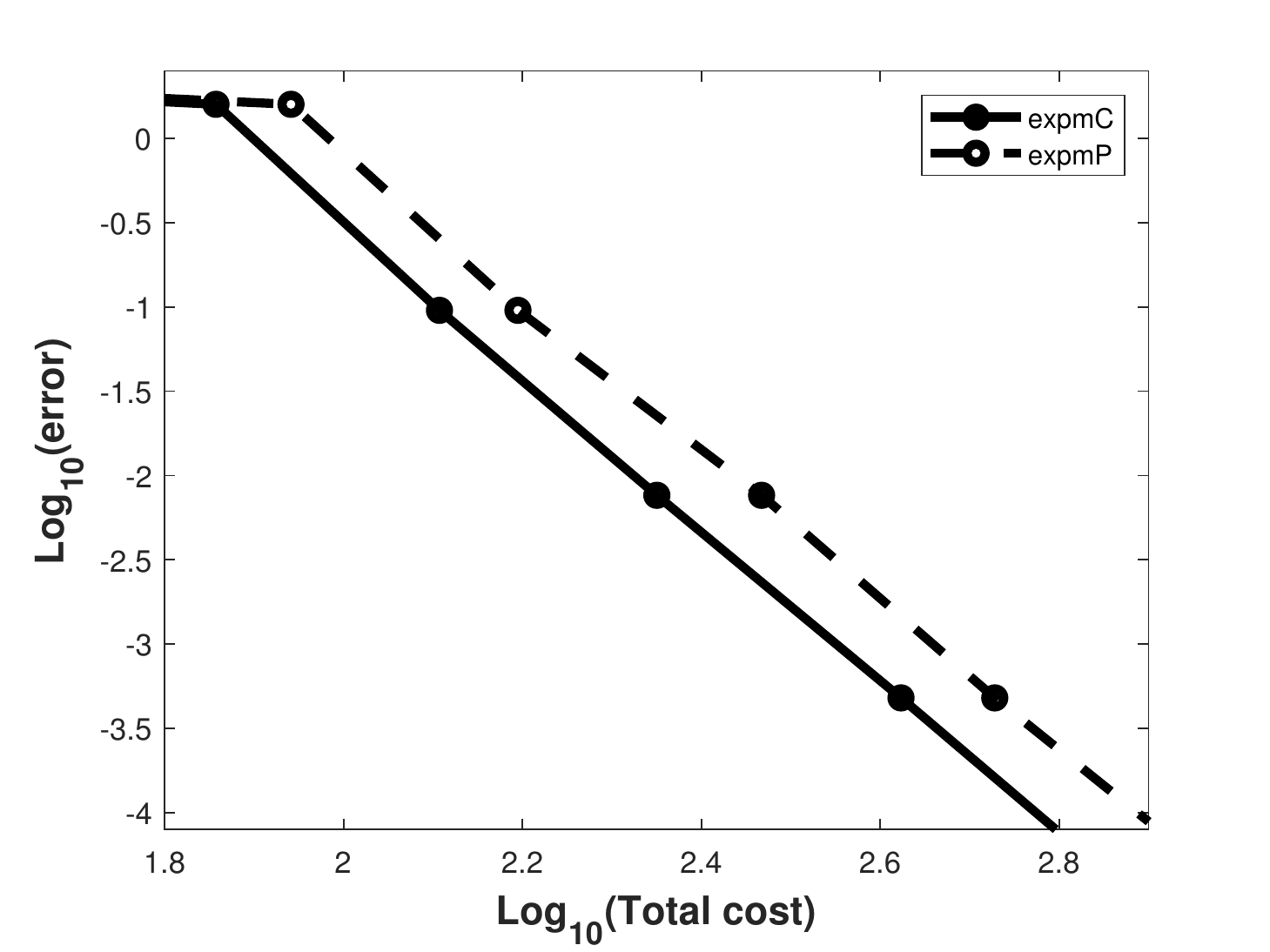}
\caption{\label{ex1err_cost}\small 2-norm error in the unitary matrix evolution at the final time versus the cost (measured as the number of matrix-matrix products required to compute the exponentials at each step) for Example 1: (top) results for the second order exponential midpoint rule and (bottom) results for the fourth order commutator-free exponential Magnus integrator.
}
\end{figure}

In this example we illustrate the performance of the new algorithm as applied to two different exponential integrators: (i) the well-known 2nd-order exponential
midpoint rule
\[
   U_{n+1} = \exp\left( -i\T H(t_n+\frac{\T}2)\right), \qquad n=0,1,2,\ldots,M-1
\]
and (ii) the 4th-order commutator-free Magnus integrator given by 
\[
   U_{n+1} = \exp\Big( -i\T (\beta H_1+\alpha H_2)\Big) \exp\Big( -i\T (\alpha H_1+\beta H_2)\Big), 
\]
where $H_i=H(t_n+c_i \T), i=1,2,$ and 
\[
 c_1=\frac12-\frac{\sqrt{3}}6, \qquad c_2=\frac12+\frac{\sqrt{3}}6, \qquad 
 \alpha=\frac14+\frac{\sqrt{3}}6,  \qquad \beta=\frac14-\frac{\sqrt{3}}6.
\] 
(See \cite{alvermann11hoc,blanes17hoc,blanes06fas} for more details of this scheme as for other higher order methods of the same class). 
The exponential matrix is computed in all cases with Pad\'e approximants and the new algorithm based on Chebyshev polynomials. Since the 
results obtained with
Taylor polynomials lie in between both of them, they are not shown in the figures for clarity). We compute
\[
  \tilde U_h = U_M \, U_{M-1} \cdots U_2 U_1,
\]
and measure the 2-norm of the error, $\|\tilde U_h-U_{\text{ref}}(t_f,t_0))\|$ for different values of $\T$. The total cost is taken as the sum of the number of the matrix-matrix products that are required for the calculation of 
$U_1 ,  \cdots,  U_M$,
and we depict the error as a function of this total number
of matrix-matrix product evaluated by each procedure. 
Figure~\ref{ex1err_cost} shows the corresponding results obtained by new procedure based on Chebyshev (expmC) and Pad\'e approximants (expmP)
for the exponential mid-point rule (top) and the 4th-order commutator-free Magnus integrator (bottom). 
We see that the relative saving in the computational cost is similar in both cases but the improvement in the accuracy increases with the order of the method.

Notice that the accuracy improves when the time step $\tau$ decreases so, the number of exponentials increases, but the cost to compute each exponential can decrease because $\|\tau H_k\|$ takes smaller values. The slope of the curves is then higher than expected from the order of the numerical integrator used.



\paragraph{Example 2: The Walker--Preston model.} 

This constitutes a standard model for a diatomic molecule in a strong laser field \cite{walker77qvc1}. The system is described by the one-dimensional Schr\"odinger equation (in units such that $\hbar=1$)
\begin{equation}\label{SchrWP}
 i \frac{\partial}{\partial t} \psi (x,t) = \left(
     -\frac{1}{2\mu} \frac{\partial^2}{\partial x^2} +
     V(x) + f(t) x \right) \psi (x,t),
\end{equation}
with $\psi(x,0)=\psi_0(x)$. Here $V(x)=D \left(1-\e^{-\alpha x}
\right)^2 $  is the Morse potential  and $f(t)x=A\cos(\omega (t))x$
accounts for the laser field. As an initial condition, we take the ground state of the Morse potential

\begin{equation}\label{incgsmor}
\psi_0(x)=\sigma \, \exp \left(-(\gamma-\frac{1}{2})\alpha x \right) \, \exp(-\gamma \, \e^{-\alpha x}),
\end{equation}
where $\gamma=2D/\omega_0$, $\omega_0=\alpha \sqrt{2D/\mu} $, and $\sigma$ is a normalizing constant.

We define the wave function $\psi$ in a certain domain $x\in[x_0,x_N]$ that is subdivided into $N$ parts of length $\Delta x=(x_N-x_0)/N$ with $x_i=x_0+i\Delta x$, and then the vector $u(t) \in \mathbb{C}^N$ with components $u_i=(\Delta x)^{1/2}\psi(x_{i-1},t), \ i=1,\ldots,N$,
is formed.

If second-order central differences are applied to discretize the equation in space and periodic boundary conditions are considered, one ends up with
the differential equation 
\[
  i \frac{du}{dt} = H(t) u = (T + B(t)) u, \qquad u(0) = u_0
\]
with
\[
  T = \frac{N^2}{2\mu(x_N-x_0)^2} \left( \begin{array}{rrrrr}
   2 & -1 & & & 1 \\
   -1 & 2 & -1 & &  \\
   & & \ddots & & \\
    & & -1 & 2 & -1 \\
   1 & & & -1 & 2
   \end{array}  \right) 
\]
and $B(t) = \mbox{diag} \big(V(x_1) + f(t) x_1, \ldots, V(x_N) + f(t) x_N \big)$. Notice that $H$ is a real symmetric matrix, $H^T=H\in\R^{N\times N}$, so that
\begin{equation} \label{eq:unitary}
   \exp\left( -i\tau H(t_n+\tau/2)\right) 
	= \cos \left( \tau H(t_n+ \tau/2) \right) -i \sin \left( \tau H(t_n+\tau/2) \right).
\end{equation}
Moreover, we can take
\[
  \Emin =  \min_{1 \le j \le N} B(t)_{jj},  \qquad   
	\Emax =  \frac{2N^2}{\mu(x_N-x_0)^2} + \max_{1 \le j \le N} B(t)_{jj}
\]
and so we shift the original matrix according with Eq. (\ref{td.2Cheb}). For our experiments we take
$x \in [-0.8, 4.32] $, the interval is subdivided into $N=64$
parts of length $\Delta x=0.08$, and the parameters are chosen as follows (in atomic units): $\mu= 1745$, $ D= 0.2251$ and $\alpha= 1.1741$
(corresponding to the HF molecule). Concerning the interaction with the laser field, we take $A= 0.011025$ and the laser frequency $\omega= 0.01787$. 

\begin{figure}[!ht] 
\centering
  \includegraphics[width=.85\textwidth]{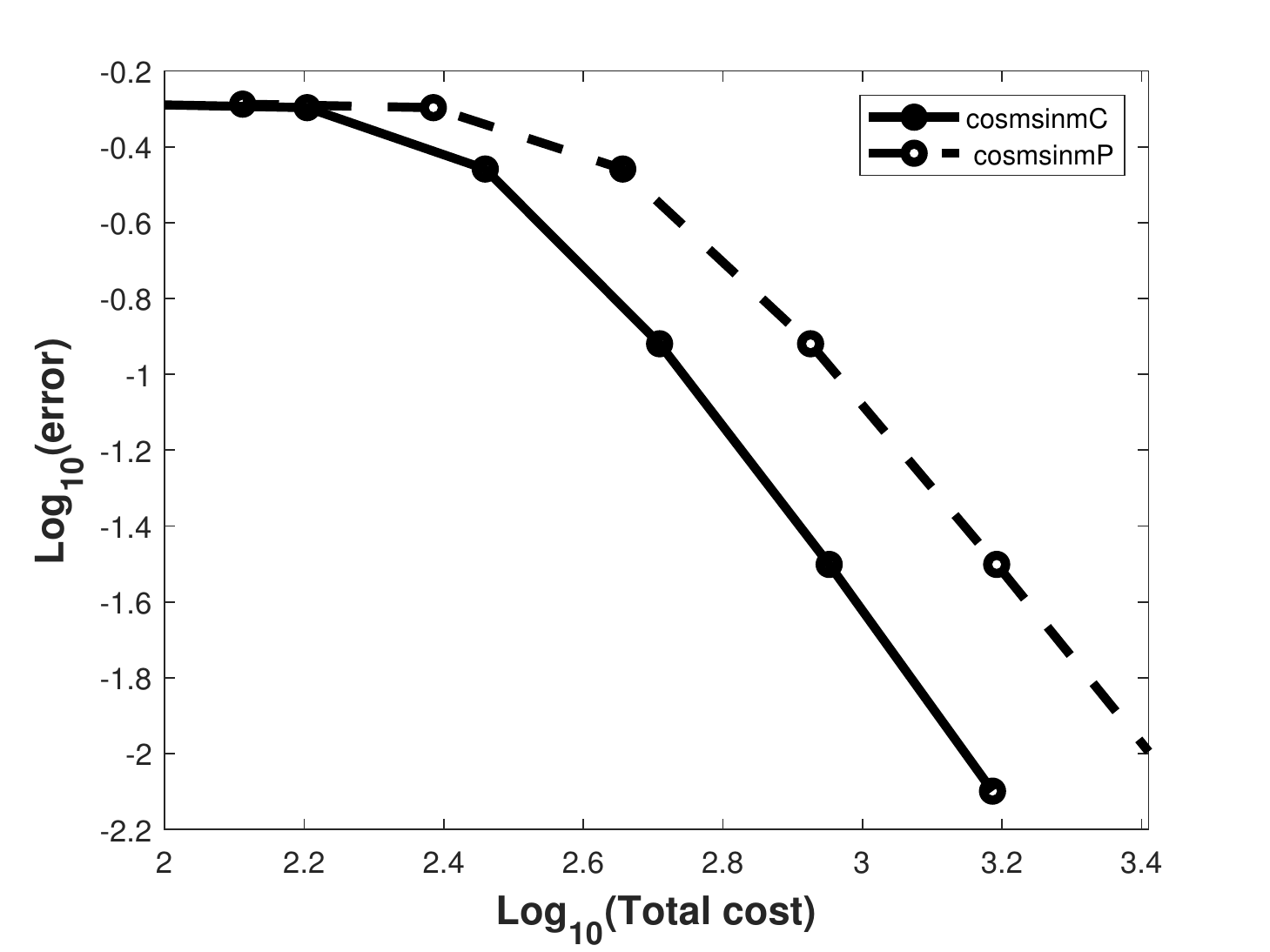}
  \includegraphics[width=.85\textwidth]{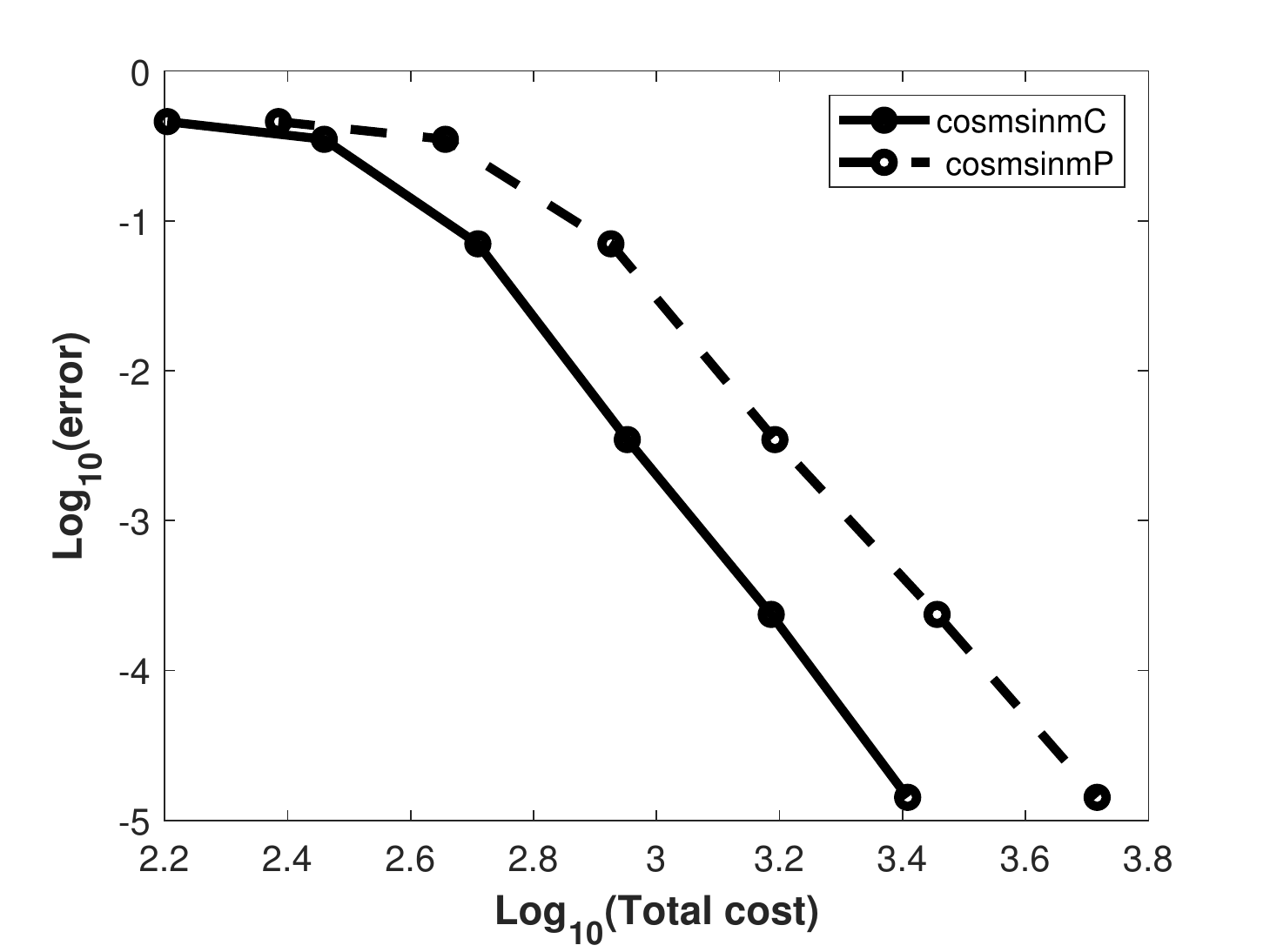}
\caption{\label{ex2err_cost}\small 2-norm error in the unitary matrix evolution at the final time versus the cost (measured as the number of matrix-matrix products required to compute the exponentials at each step) for the real symmetric matrix $H$ of Example 2: (top) results for the second order exponential mid-point rule and (bottom) results for the fourth order commutator-free exponential Magnus integrator.
}
\end{figure}

As before, to check the performance of the different procedures, we compute the 2-norm error in the evolution matrix solution at the final time
$t_f = \frac{2\pi}{\omega}$. To do that, we compare with a reference solution computed with high accuracy. The total cost of each procedure
is measured as the total number of matrix-matrix products required to approximate the matrix cosine and sine for the total integration interval.  In this way
we get Figure~\ref{ex2err_cost}, where the results achieved by Chebyshev approximations (denoted by `cosmsinmC') and Pad\'e approximants
(`cosmsinmP', obtained with the algorithm of \cite{almohy15naf}) are collected. 
The top diagram corresponds to the 2nd-order exponential mid-point rule and the bottom  graph is obtained with the 4th-order commutator-free Magnus
integrator.
Here again, the new algorithm based on Chebyshev polynomials leads to more accurate results with a reduced computational cost. 




\section{Conclusions and future work} \label{sec.5}

We have presented an algorithm to approximate the exponential of skew-Her\-mi\-tian matrices  based on an improved computation of Chebyshev polynomials of matrices and the corresponding error analysis. For problems of the form $\exp(-iA)$, when $A$ is a real and symmetric matrix, an improved version is presented that computes the sine and cosine of $A$ 
with a reduced number of products of real and symmetric matrices. In both cases, the new procedures turn out to be more 
efficient than schemes based on rational Pad\'e approximants or Taylor polynomials for all
tolerances and time interval lengths.

The Chebyshev methods presented in this paper can be further improved along different lines that will be explored in our future work:

\begin{itemize}
 \item As we have seen, with only three products it is possible to evaluate most polynomials to order eight (this is, in fact, 
 the highest degree one can reach with three products). With four products one can build polynomials of degree sixteen, but there are not enough free parameters to obtain the Taylor and Chebyshev polynomials approximating the exponential.  For this reason, we have limited ourselves
here to polynomials of order twelve, which can be obtained with four products. On the other hand, in \cite{sastre19btc} a polynomial of degree
16 is presented in terms of only 4 products that coincides with the Taylor expansion up to order 15 (this method is denoted in Table~\ref{tab.thetaTaylor} as $m=15+$). In this way, with the computational cost as the method of degree $m=12$, it 
provides a larger value for $\theta$ that is even slightly larger that the value of the Chebyshev polynomial of degree 12. The same procedure
can of course be carried out with Chebyshev polynomials: one could construct a polynomial of degree 16 that coincides with the Chebyshev polynomial up to degree 15 and analyze whether this new polynomial has a larger value of $\theta$. Notice that the procedure is largely
similar to the search of polynomials $S_{m,\theta}^C$ that coincide with $s_{m,\theta}^C$. With five products it is also possible to build a polynomial of degree 24 that approximates the Chebyshev polynomial up to order 21, and we expect an improvement with respect to the result obtained for $m=18$, in the same way as in \cite{sastre19btc} for the Taylor polynomial.

\item One could also build a new set of methods aimed to be used with different accuracies, and in particular in single precision.
 From the error bound formulas for the chosen values of $m$, the new values for $\theta$ have to be obtained and then the corresponding Chebyshev polynomials of degree $m$ have to be obtained that will be then computed with a reduced number of products.

\item When lower accuracies are desired then the preservation of unitarity is also lost to such accuracy. It is well known that diagonal Pad\'e methods preserve unitarity unconditionally and one can look for similar rational Chebyshev approximations to analyze the preservation of unitarity as well as to reduce the cost of these schemes. Rational Chebyshev approximations have been successfully used in \cite{sidje98eas} to compute the action of the exponential of skew-Hermitian matrices on vectors.

\item Finally, there are of course a number of efficient procedures for the diagonalization of Hermitian or skew-Hermitian matrices that might be also employed for
evaluating the matrix exponentials required for the application to exponential integrators to certain classes of differential equations. In that case the norm
of the matrices involved is usually quite small (since they involve the step size of the integrator) and thus our algorithms are particularly well suited for this
purpose. In any case, a future line of research consists in determining precisely under which circumstances related with the size and norm of the matrix 
the algorithms presented here are competitive with other procedures based on direct diagonalization.
\end{itemize}

\subsection*{Acknowledgements}

SB and FC have been supported by Ministerio de Ciencia e Innovaci\'on (Spain) through project  PID2019-104927GB-C21 (AEI/FEDER, UE). The work of MS has been funded by the Scientific and Technological Research Council of Turkey (TUBITAK) with grant number 1059B191802292. 
SB and FC would like to thank the Isaac Newton Institute for Mathematical Sciences for support and hospitality during the programme ``Geometry, compatibility and structure preservation in computational differential equations'', when work on this paper was undertaken. This work was been additionally supported by EPSRC grant number EP/R014604/1. The authors wish to thank the referee for his/her detailed list of comments and
suggestions which were most helpful to improve the presentation of the paper.

\bibliographystyle{siam}

\end{document}